\documentclass[a4paper, 12pt, leqno]{article}%{amsart}

\usepackage{latexsym, amsmath, amsthm, amssymb, graphicx, graphics, epsfig, bm, fix-cm}
\usepackage{tikz-cd}

\usepackage[T1]{fontenc}

\newcommand{\comment}[1]{}

\usepackage{fullpage}

\usepackage{manfnt}
\usepackage[sans]{dsfont}
\usepackage{palatino}
\usepackage{eulervm}
\usepackage{euscript}

\usepackage{dsfont} %{bbm}

\newlength{\otstup} 
\newcommand{\sub}[1]{\vspace{\otstup}\textbf{#1}\hspace*{0.5em}}
\newcommand{\upskip}{\vspace{-\otstup}}
\newcommand{\uppskip}{\vspace{-2\otstup}}

\newtheorem{Thm}{Theorem}
\usepackage{xpatch}
\xpatchcmd\swappedhead{~}{.~}{}{}

\newtheorem{thm}{Theorem}[section]
\newtheorem*{thm*}{Theorem}
\newtheorem{lem}[thm]{Lemma}
\newtheorem*{lem*}{Lemma}

\newtheorem*{cor*}{Corollary}
\newtheorem{prop}[thm]{Proposition}
\newtheorem*{prop*}{Proposition}
\newtheorem*{claim*}{Claim}

\theoremstyle{definition}

\newtheorem{rem}[thm]{Remark}
\newtheorem*{rem*}{Remark}

\newtheorem{exm}[thm]{Example}

\numberwithin{equation}{section}

\renewcommand{\proof}{\vspace{-8pt}\noindent\textit{\textbf{Proof. }}}

\renewcommand{\endproof}{$\square$}

\newcommand{\txt}[1]{\;\;\text{#1}\;\;}
\newcommand{\q}[1]{{``#1''}}
\newcommand{\br}[1]{\left(#1\right)}

\newcommand{\bra}[1]{\left\langle #1 \right\rangle}

\renewcommand{\leq}{\leqslant}
\renewcommand{\geq}{\geqslant}

\newcommand{\tor}[1]{\stackrel{{#1}\hspace{0.2em}}{\longrightarrow}}
\newcommand{\isor}{\tor{\approx}}

\newcommand{\hookto}{\hookrightarrow}

\newcommand{\set}[1]{\left\{#1\right\}}
\newcommand{\sett}[2]{\left\{#1 \mid #2\right\}}

\newcommand{\norm}[1]{\lVert #1 \rVert}

\newcommand{\inv}{^{-1}}
\newcommand{\restr}[2]{\left. #1 \right|_{#2}}
\newcommand{\eps}{\varepsilon}
\renewcommand{\phi}{\varphi}

\newcommand{\matr}[1]{
\begin{pmatrix}
#1
\end{pmatrix}}

\newcommand{\smatr}[1]{\br{ \begin{smallmatrix}#1\end{smallmatrix} } }

\newcommand{\R}{\mathbb R}
\newcommand{\N}{\mathbb N}

\newcommand{\C}{\mathbb C}

\newcommand{\PP}{\mathbb P}

\DeclareMathOperator{\Res}{Res}

\DeclareMathOperator{\Ran}{Ran}

\DeclareMathOperator{\codim}{codim}

\DeclareMathOperator{\Id}{Id}

\DeclareMathOperator{\dom}{dom}

\newcommand{\sa}{^{\mathrm{sa}}}

\newcommand{\Gr}{\mathcal Gr}

\newcommand{\U}{\mathcal U}
\newcommand{\sU}{{^s{\mathcal U}}}
\newcommand{\sP}{{^s{\mathcal P}}}
\newcommand{\sB}{{^s{\mathcal B}}}

\newcommand{\xX}{x\in X}

\newcommand{\Reg}{\mathcal R}
\newcommand{\Rsa}{\Reg\sa}
\newcommand{\RF}{\Reg_F}
\newcommand{\RK}{\Reg_K}

\newcommand{\rR}{{^r\Reg}}
\newcommand{\gR}{{^g\Reg}}

\newcommand{\D}{\mathcal D}

\newcommand{\B}{\mathcal B}

\newcommand{\BF}{\B_F}

\newcommand{\Proj}{\mathcal P}
\newcommand{\PK}{\Proj_{K}}

\newcommand{\st}{^{\star}}
\newcommand{\Po}{\Proj\st}

\newcommand{\Ro}{\Reg\st}

\DeclareMathOperator{\one}{\mathds{1}}

\newcommand{\onelp}{\one_{[\lambda,+\infty)}}

\newcommand{\K}{\mathcal K}
\renewcommand{\H}{\mathcal H}

\newcommand{\M}{\mathcal M}
\newcommand{\T}{\mathcal T}
\newcommand{\F}{\mathcal F}
\newcommand{\E}{\mathcal E}
\newcommand{\A}{\mathcal A}

\newcommand{\I}{\mathcal I}
\renewcommand{\S}{\mathcal S}

\DeclareMathOperator{\ind}{ind}
\DeclareMathOperator{\Ind}{Ind}

\newcommand{\f}{\bm{\chi}}

\newcommand{\p}{\mathbf p}

\newcommand{\pinf}{p_{\infty}}
\newcommand{\w}{\mathbf w}

\newcommand{\s}{\mathbf s}

\newcommand{\ppi}{\bm{\pi}}

\title{From graph to Riesz continuity} 
\author{Marina\,Prokhorova}
\date{}

\begin{document}

\maketitle

\footnotetext{\hspace*{-1.8em}Department of Mathematics, University of Haifa (Israel) \\
Department of Mathematics, Technion -- Israel Institute of Technology \vspace*{0.5em}}

\renewcommand{\baselinestretch}{1.00}
\selectfont

\upskip
\begin{abstract}
\parindent=0cm
\setlength{\parskip}{6pt plus 0pt minus 0pt}
\noindent
We show that every graph continuous family of unbounded operators in a Hilbert space becomes Riesz continuous
after one-sided multiplication by an appropriate family of unitary operators.
This result provides a simple definition of the index for graph continuous families of Fredholm operators,
and we show that for such families this index coincides with the index defined by N.~Ivanov in \cite{I2}.
This result also has two corollaries for operators with compact resolvents: 
(1) the identity map between the space of such operators with the Riesz topology
and the space of such operators with the graph topology is a homotopy equivalence; 
(2) every graph conti\-nuous family of such operators acting between fibers of Hilbert bundles 
becomes Riesz continuous in appropriate trivializations of the bundles.

For self-adjoint operators, multiplication by unitary operators should be replaced by conjugation.
In general, a graph continuous family of self-adjoint operators  with com\-pact resolvents
cannot be made Riesz continuous by an appropriate conjugation.
We obtain a partial analogue of the \q{trivialization} result above for self-adjoint operators
and describe obstructions to existence of such a trivialization in the general case.
This motivates the notion of a polarization of a Hilbert bundle,
and we prove a similar result for polarizations.
These results  are closely related to the recent work \cite{I2} of Ivanov 
and provide alternative proofs for some of his results.
We then show that, under a minor assumption on the space of parameters
and for operators which are neither essentially positive nor essentially negative, 
there is always a trivialization making the family Riesz continuous.
\end{abstract}

\section{Introduction}\label{sec:intro}

Let $H$ and $H'$ be separable complex Hilbert spaces of infinite dimension.
We denote by $\B(H,H')$ the space of bounded linear operators $H\to H'$ with the norm topology;
by $\K(H,H')$ the subspace of $\B(H,H')$ consisting of compact operators, 
by $\U(H)$ and $\Proj(H)$ the subspaces of $\B(H)=\B(H,H)$ 
consisting of unitary operators and projections respectively 
(by projections we always mean orthogonal projections). 

\sub{Regular operators.}
An unbounded operator $A$ from $H$ to $H'$ is a linear operator 
defined on a subspace $\dom(A)$ of $H$ and taking values in $H'$.
Such an operator $A$ is called closed if its graph is closed in $H\oplus H'$, 
and densely defined if its domain $\dom(A)$ is dense in $H$.
It is called \emph{regular} if it is closed and densely defined.
Let $\Reg(H,H')$ denote the set of all regular operators $H\to H'$
and $\Reg\sa(H)\subset \Reg(H) = \Reg(H,H)$ denote the subset of self-adjoint operators.

\sub{Two topologies on regular operators.}
The set $\Reg(H,H')$ of regular operators has two important natural topologies.
The \emph{Riesz topology} on $\Reg(H,H')$ is induced by the so called \emph{bounded transform map} 
\[ \f\colon\Reg(H,H')\hookto\B(H,H'), \quad \f(A) = A(1+A^*A)^{-1/2} \]
from the norm topology on the space of bounded operators.
The \emph{graph topology} on $\Reg(H,H')$ is induced by the inclusion 
\[ \p\colon\Reg(H,H')\hookto\Proj(H\oplus H'), \] 
from the norm topology on the space $\Proj(H\oplus H')$ of projections, 
where $\p$ is the map taking a regular operator to the orthogonal projection onto its graph.
Let $\rR(H,H')$, resp. $\gR(H,H')$ denote the space of regular operators $H\to H'$ equipped with the Riesz, resp. graph topology.

The Riesz topology is strictly finer than the graph topology.
In other words, the identity map $\rR\to\gR$ is continuous, while the identity map $\gR\to\rR$ is not.
On the subset $\B(H,H')$ of bounded operators these two topologies coincide.

\sub{From graph to Riesz continuity.} 
Using the bounded transform, one can deal with Riesz continuous families of operators
in essentially the same manner as with norm continuous families of bounded operators.
Because of that, Riesz continuous families are much easier to deal with than mere graph continuous families.
Our first result shows that one can turn a graph continuous family to a Riesz continuous one
by a \q{change of coordinates}\footnote{The interpretation of our construction as a change of coordinates was suggested by the referee.}.

\begin{Thm}\label{Thm:Phi}
	There exist a map $v\colon\Reg(H,H')\to\U(H)$, which is both Riesz-to-norm and graph-to-strong continuous,
	such that the induced map 
	\[ \Phi\colon\Reg(H,H')\to\Reg(H,H'), \quad A\mapsto A\cdot v(A) \] 
	is graph-to-Riesz continuous. 
\end{Thm}

\upskip
Passing to adjoint operators, one obtains a similar statement for the left multiplication.
Namely, there is a map $v'\colon\Reg(H,H')\to\U(H')$, which is both Riesz-to-norm and graph-to-strong continuous,
such that the induced map $A\mapsto v'(A)\cdot A$ is graph-to-Riesz continuous. 
\medskip

We will now present three applications of Theorem \ref{Thm:Phi}.

\sub{$K^0$ index for graph continuous families.}
Let $\A\colon X\to\gR$ be a graph continuous map
and $\Phi\colon\gR\to\rR$ be a map provided by Theorem \ref{Thm:Phi}.
Then the composition 
\[ \alpha=\f\circ\Phi\circ\A\colon X\to\B(H) \] is norm continuous.
Suppose that all $\A_x$ are Fredholm. 
Then the operators $\alpha_x$ are also Fredholm, and one can define the index of $\A$ 
as the classical index of $\alpha\colon X\to\B_F$
(here $\B_F$ denotes the space of bounded Fredholm operators):
\begin{equation*}
	\ind\A = \ind\alpha = [\alpha] \in [X,\BF] = K^0(X).
\end{equation*}
This definition of the index is justified by Theorem \ref{Thm:RKF} below.
We discuss the index of graph continuous families in the last part of Section \ref{sec:appl1}.

Searching for such a definition of the index was part of the motivation behind Theorem~\ref{Thm:Phi}.
An alternative approach to the index of graph continuous families of Fredholm operators
is provided by the previous paper \cite{Pr21} of the author.
Another approach, for compact base spaces and in a more general framework of Hilbert modules over C*-algebras,
is provided by results of M.~Joachim \cite{Jo}.
Recently N.~Ivanov developed a new approach to the index of families \cite{I1, I2}, 
which works under much weaker continuity assumptions than graph continuity. 
We show in the end of Section \ref{sec:appl1} that the definition of the index given above 
coincides with Ivanov's index for graph continuous families of Fredholm operators.

\sub{Homotopy equivalences.}
A regular operator $A$ is said to have compact resolvents\footnote{If the resolvent set 
$\Res(A)$ of $A$ is non-empty (for example, if $A$ is self-adjoint), 
then this definition agrees with the usual sense of the words \q{compact resolvent}:
$(A-\lambda)\inv$ is compact for some (and then every) $\lambda\in\Res(A)$.
However, a regular operator may have an empty resolvent set; 
the definition above covers such operators as well.} 
if both $(1+A^*A)\inv$ and $(1+AA^*)\inv$ are compact operators.
Let $\RF(H)$ and $\RK(H)$ be the subsets of $\Reg(H)$ consisting of Fredholm operators and 
operators with compact resolvents, respectively.
These subsets are invariant under the multiplication by unitary operators, 
so $\Phi(\RK)\subset\RK$ and $\Phi(\RF)\subset\RF$.

It was shown in the previous paper of the author that the identity maps 
\begin{equation*}
	\rR_K\to\gR_K \txt{and} \rR_F\to\gR_F 
\end{equation*}
are homotopy equivalences. See \cite[Theorem A]{Pr21}.
However, the proof in \cite{Pr21} provides no intuition on how a homotopy inverse map might look like.
The following theorem reveals a possible form of homotopy inverse maps and
provides an alternative proof of homotopy equivalence of the identity map $\rR_K\to\gR_K$.

\begin{Thm}\label{Thm:RKF}
Let $\Phi\colon\gR\to\rR$ be a map provided by Theorem \ref{Thm:Phi}.
Then the restriction of $\Phi$ to the subspace of operators with compact resolvents 
is homotopy inverse to the identity map $\rR_K\to\gR_K$.
Similarly, the restriction of $\Phi$ to the subspace of Fredholm operators 
is homotopy inverse to the identity map $\rR_F\to\gR_F$.
\end{Thm}

\upskip
\sub{Hilbert bundles.}
Let $\H$ and $\H'$ be locally trivial Hilbert bundles over a topological space $X$ with a fiber $H$.
(In the following we will always suppose that all Hilbert bundles are locally trivial.)
The structure group of a general Hilbert bundle is $\sU(H)$, 
the unitary group of $H$ with the strong operator topology.
This group is contractible \cite[Lemma 3]{DD},
so every such bundle over a paracompact base space $X$ is trivial, 
that is, isomorphic to the trivial Hilbert bundle $X\times H\to X$. 
More generally, every numerable Hilbert bundle is trivial.
However, automorphisms of the trivial Hilbert bundle $X\times H\to X$ are given by continuous maps $X\to\sU(H)$
and preserve neither norm nor Riesz continuity of operator families, 
since the action of $\sU(H)$ on $\B(H)$ by conjugation, $(u,a)\mapsto uau\inv$, is not continuous.
Therefore, the notions of a norm continuous family of bounded operators 
or a Riesz/graph continuous family of regular operators
are not defined for a general Hilbert bundle.

Nevertheless, the actions of $\sU(H)$ and $\sU(H')$ on the space $\K(H,H')$ of compact operators 
are continuous, and thus the notion of a norm continuous family of compact operators 
$k_x\colon\H_x\to\H'_x$ is well defined.
By the same reason, the notion of a \textit{graph continuous} family of regular operators 
$\A_x\colon\H_x\to\H'_x$ with \textit{compact resolvents} is well defined (Proposition \ref{prop:sURK}).
Applying Theorem \ref{Thm:Phi} to such a family, we obtain the following result.

\begin{Thm}\label{Thm:RK-H}
Let $\H$ and $\H'$ be numerable Hilbert bundles over a topological space $X$.
Let $\A$ be a graph continuous family of regular operators $\A_x\colon\H_x\to\H'_x$ with compact resolvents.
Suppose that a trivialization of one of these Hilbert bundles, 
$\H\cong X\times H$ or $\H'\cong X\times H'$, is fixed. 
Then there is a trivialization of the second bundle taking $\A$ 
to a Riesz continuous map $X\to\RK(H,H')$.	
\end{Thm}

\upskip
\sub{Self-adjoint case.} 
For self-adjoint operators, multiplication $A\mapsto Au$ or $A\mapsto uA$ 
should be replaced by the conjugation $A\mapsto uAu\inv$.
Instead of independent trivializations of two Hilbert bundles $\H$ and $\H'$, as in Theorem \ref{Thm:RK-H},
one should look for a one trivialization of a given Hilbert bundle $\H$.
This narrows down available options drastically.
In particular, Theorem \ref{Thm:Phi} has no direct analogue for self-adjoint operators, even invertible ones.

Let $\A$ be a family of regular operators acting on fibers of a locally trivial Hilbert bundle $\H$ over $X$.
Following Ivanov \cite{I3}, we say that a trivialization $\tau\colon\restr{\H}{Y}\to Y\times H$ of $\H$ 
over a subset $Y\subset X$ is \emph{adapted}\footnote{\q{Fully adapted} in terms of Ivanov} to $\A$ 
if it takes $\A$ to a Riesz continuous map $\tau_*\A\colon Y\to\Reg(H)$.
In particular, $\tau$ is adapted to a family $\A$ of \emph{bounded} operators 
if it takes $\A$ to a \emph{norm} continuous map $\tau_*\A\colon Y\to\B(H)$.

Spectral projections of a graph continuous family $\A$ of self-adjoint operators with compact resolvents
give rise to what we call a compatible family of local bundle projections.
A trivialization is adap\-ted to $\A$ if and only if it is adapted to these projections.
Our main results in the self-adjoint case, Theorems \ref{Thm:Gr} and \ref{Thm:Gr-fd},
are concerned with such families of projections,
and we deduce results concerning self-adjoint operators from these two theorems.

\sub{Polarizations.} 
We say that $p\colon\restr{\H}{Y}\to\restr{\H}{Y}$ is a \emph{bundle projection of $\H$ over $Y\subset X$}
if $p$ is a continuous bundle map such that each $p_x$ is a projection.
We say that two bundle projections, $p$ over $Y$ and $p'$ over $Y'$, are \emph{compatible} 
(resp. \emph{finitely compatible})
if the difference $x\mapsto p_x-p'_x$ is a norm continuous family of compact (resp. finite rank) 
operators over $Y\cap Y'$.

Locally, a polarization is defined by a bundle projection $p$ over an open subset $Y\subset X$ 
such that each $p_x$ has infinite-dimensional range and kernel.
Globally, a polarization is defined by an open covering $(X_i)$ of $X$ together with a collection 
of pairwise compatible local polarizations $p_i$ over $X_i$.
Following Ivanov \cite{I2}, we call such a collection $(X_i,p_i)$ a \textit{polarization atlas}
and define a \textit{polarization} as an equivalence class of polarization atlases
(as usual, two atlases are equivalent if their union is also an atlas).
A \emph{finite polarization} is defined in the same way, but with \q{compatible} replaced by \q{finitely compatible}.
See Section \ref{sec:polar} for details.

Let $\Pi$ be a polarization of $\H$ given by an atlas $(X_i,p_i)$.
We say that a (local or global) trivialization of $\H$ is \emph{adapted} to $\Pi$
if it is adapted to each $p_i$ (that is, takes each $p_i$ to a norm continuous family).
We say that a polarization is \emph{tame} if it admits a local adapted trivialization 
over a neighborhood of each point $x\in X$.

The notion of a polarization introduced by Ivanov in the appendix to \cite{I2}
is essentially equivalent to a tame finite polarization in our terminology.

\begin{Thm}\label{Thm:Gr}
Let $\Pi$ be a tame polarization of a locally trivial Hilbert bundle $\H$ over a paracom\-pact space $X$.
Then $\Pi$ admits a global adapted trivialization.
Moreover, every adapted trivialization of $\H$ over a closed subset $Y$ of $X$ 
can be extended to a global adapted trivialization.
In particular, every two adapted trivializations are homotopic (as trivializations adapted to $\Pi$).
\end{Thm}

The next theorem shows that a modest assumption on the base space is enough 
for existing of an adapted trivialization.

Let us say that a paracompact space $X$ has \emph{locally finite covering dimension} 
if each point of $X$ has a neighborhood of finite covering dimension
(we do not require this neighborhood to be open or small).
In particular, every space of finite covering dimension has locally finite covering dimension.

\begin{Thm}\label{Thm:Gr-fd}
Let $\H$ be a locally trivial Hilbert bundle over a paracompact space $X$ of locally finite covering dimension.
Then every polarization of $\H$ is tame and thus admits a global adapted trivialization.
\end{Thm}

\upskip
\sub{Self-adjoint operators and adapted trivializations.} 
We want to know whether a family $\A$ of self-adjoint operators $\A_x\colon\H_x\to\H_x$
admits a global adapted trivialization.
Clearly, for that to be true, $\A$ has to admit local adapted trivializations
(by that we mean that for every $x\in X$ there is an adapted trivialization over a neighborhood of $x$). 
Theorem \ref{Thm:Gr} implies the following result, 
which shows that for essentially unitary operators this local condition is also sufficient.
The author does not know whether this is true for arbitrary bounded self-adjoint operators, even if they are Fredholm.

\begin{Thm}\label{Thm:Beusa}
  Let $\H$ be a locally trivial Hilbert bundle over a paracompact space $X$,
	$\A$ be a family of self-adjoint operators $\A_x\colon\H_x\to\H_x$,
	and either all operators $\A_x$ are essentially unitary or they all have compact resolvents. 
	Suppose that $\A$ admits local adapted trivializations. 
	Then $\A$ admits a global adapted trivialization.
	Moreover, every adapted trivialization of $\H$ over a closed subset $Y$ of $X$ 
	can be extended to a global adapted trivialization.
	In particular, every two such trivializations are homotopic (in the space of trivializations adapted to $\A$). 
\end{Thm}

The main case of Theorem \ref{Thm:Beusa} is the case of families consisting of operators 
which are neither essentially positive nor essentially negative. 
In this case, a similar trivialization theorem 
was nearly simultaneously proved by different methods and under much weaker continuity assumptions by Ivanov. 
See \cite{I2}, Theorems 4.5--4.7. 
It can be easily shown that his Theorems 4.6--4.7 imply Theorem \ref{Thm:Beusa} 
for families consisting of operators which are neither essentially positive nor essentially negative.
Ivanov mentioned in the introduction to the first version of \cite{I2} 
that there is a notion of a polarization of a Hilbert bundle standing behind some proofs in \cite{I2}. 
This led the author to introducing her version of a polarization of a Hilbert bundle 
and proving an analogue of Ivanov's trivialization theorems 4.6--4.7 in terms of polarizations (Theorem \ref{Thm:Gr}).
Later, Ivanov wrote down the details of his notion of a polarization 
in the appendix to the current version of \cite{I2} 
and rephrased the proofs of his Theorems 4.5--4.7 in terms of polarizations. 
See \cite{I2}, Theorems A.2--A.4.

\sub{First local obstruction.}
Let $\Rsa_K(H)$ denote the set of self-adjoint operators with compact resolvents.
The first local obstruction to existence of an adapted trivialization 
for a graph continuous family $\A$ of operators $\A_x\in\Rsa_K(\H_x)$
is caused by the fact that $\Rsa_K$ is path connected in the graph topology,
but has three connected components in the Riesz topology.
Two points from different components of $\rR_K\sa$
(for example, an essentially positive and an essentially negative operator)
may be connected by a graph continuous path in $\RK\sa(H)$,
but no trivialization of the trivial Hilbert bundle $[0,1]\times H$ makes this path Riesz continuous.

\sub{Second local obstruction.} 
When all the operators $\A_x$ are contained in a fixed connected component of $\rR_K\sa$,
there is another kind of local obstruction which we now describe.
For a fixed $x\in X$, choose $\lambda\in\R$ in the resolvent set of $\A_{x}$.
Then $\lambda$ also lies in the resolvent set of $\A_y$ for $y$ close to $x$, 
and the spectral projections $p_y = \onelp(\A_y)$ determine a continuous map (a bundle projection) 
$p\colon\restr{\H}{Y}\to\restr{\H}{Y}$ over some neighborhood $Y$ of $x$
(here $\one_S$ denotes the characteristic function of a subset $S\subset\R$).
See Proposition \ref{prop:p-g2s}.
A trivialization of $\H$ over $Y$ is adapted to $\A$ if and only if it is adapted to $p$
(Proposition \ref{prop:ad-pr}).

Suppose first that all the operators $\A_x$ are essentially negative. 
On the space of projections of a fixed finite rank the norm and strong topologies coincide \cite[Lemma~5]{DD},
so there is an obstruction to existence of local adapted trivializations near $x$
if the rank of $p_y$ is not locally constant near $x$.
In other words, there are eigenvalues of $\A_y$ going to $+\infty$ when $y\to x$. 
If this is not the case and the operators $\A_y$ are uniformly bounded from above 
in a neighborhood of each point $x\in X$,
then $\A$ is Riesz continuous in every trivialization of $\H$.
The situation with essentially positive operators is completely the same, 
but $+\infty$ and bounds from above should be replaced by  $-\infty$ and bounds from below.
A graph continuous family of essentially positive (resp, essentially negative) self-adjoint operators with compact resolvents may have no local adapted trivializations in a neighborhood of $x\in X$, 
even when the base space $X$ is an interval.
See Rellich's example \cite[Example V-4.14]{Kato}, or \cite[Example 8.2]{Pr21} for more details.

Suppose now that all $\A_x$ are contained in the third connected component $\rR_K\st(\H_x)$ of $\rR_K\sa(\H_x)$, 
which consists of operators that are neither essentially positive nor essentially negative.
In general, such a family, as well as its spectral sections, may have no local adapted trivializations. 
See Example \ref{exm:DD}.
However, Theorem \ref{Thm:Gr-fd} implies the following result.

\begin{Thm}\label{Thm:RKsa-fd}
  Let $\H$ be a locally trivial Hilbert bundle over a paracompact space $X$ of locally finite covering dimension.
	Let $\A$ be a graph continuous family of self-adjoint operators $\A_x\colon\H_x\to\H_x$, 
	with compact resolvents, which are neither essentially positive nor essentially negative.
	Then $\A$ admits a global adapted trivialization.
\end{Thm}

\upskip
Since the group $\sU(H)$ is contractible, 
Theorem \ref{Thm:RKsa-fd} implies that the identity map $\rR\st_K\to\gR\st_K$ is a weak homotopy equivalence. 
The stronger result, of a homotopy equivalence of this map, 
is proven in the previous paper of the author by a different method.
See \cite{Pr21}, Theorems B and C.	

Example \ref{exm:DD} below shows that Theorem \ref{Thm:RKsa-fd} does not hold for general base spaces, 
even compact, without additional assumptions.

The author does not know whether an analogue of Theorem \ref{Thm:RKsa-fd} 
(probably with stronger assumptions on $X$ and for maps $\A\colon X\to\gR\st_F$) 
holds for Fredholm operators.
See Remark \ref{rem:Fr-sa}.

\sub{$K^1$ index for graph continuous families.}
Let $\A$ be a family satisfying assumptions of Theorem \ref{Thm:RKsa-fd}.
Then $\A$ admits an adapted trivialization $\tau\colon\H\to H_X$, 
which leads to a norm continuous map $\alpha = \f(\tau_*\A)\colon X\to\BF\st(H)$.
The homotopy class of $\alpha$ does not depend on the choice of $\tau$.
This allows to define the index of $\A$ as the classical index of $\alpha$,
$\ind\A := \ind\alpha = [\alpha]\in [X,\BF\st] = K^1(X)$. 

Another approach to the index of graph continuous families of self-adjoint Fredholm operators
is provided by the previous paper \cite{Pr21} of the author or by results of Joachim \cite{Jo}.
An approach of Ivanov \cite{I1, I2} works even under much weaker continuity assumptions than graph continuity. 
The same reasoning as for the $K^0$ index (Section \ref{sec:appl1}) 
shows that, for a family admitting an adapted trivialization,
the index given by the formula above coincides with Ivanov's index.

\sub{Applications to elliptic operators.}
The author arrived at Theorem \ref{Thm:Beusa} in an attempt to understand a step 
in the proof of a result of Melrose and Piazza.
See \cite{MP1}, Proposition 1.
They deal with a family $\D$ of self-adjoint differential operators $\D_x$ of order $1$
parametrized by $x\in X$ such that $\D_x$ acts on the fibers of a vector bundle $\E_x$ 
over a smooth closed manifold $\M_x$.
More precisely, $\M\to X$ is a smooth fibration and $\E$ is a smooth vector bundle over $\M$.
Melrose and Piazza showed that if the $K^1(X)$-index of $\D$ vanishes, then $\D$ admits a spectral section.
In their proof, they 
claim that the $L^2$-spaces of the fibers, $\H_x = L^2(\M_x;\E_x)$, 
can be identified so that all operators act on a fixed Hilbert space $H$.
Unfortunately, they did not explain how such an identification can be obtained.
The author was not able to find a simple explanation of this step.
The structure group of the corresponding Hilbert bundle $\H$ has the topology coarser than the norm topology, 
so Kuiper's theorem cannot be used to trivialize $\H$.
Dixmier-Douady theorem, on the other hand, can be used, 
but trivializations it provides do not take $\D$ to a Riesz continuous family of operators.
At the same time, the proof of Melrose and Piazza uses essentially Riesz continuity of $\D$.

Theorem \ref{Thm:Beusa}, as well as recent results \cite[Theorems 4.5 or 4.6]{I2} 
of Ivanov, allows to overcome this difficulty by turning $\D$ to a Riesz continuous family 
(though at the cost of losing the differential and geometric structure). 
Similarly, Theorem \ref{Thm:RK-H} allows to pass to a fixed Hilbert space in the proof of 
\cite[Proposition 2]{MP2}, an odd analog of \cite[Proposition 1]{MP1}.

Another framework where our results may be useful is that of elliptic boundary value problems.
A continuous, in an appropriate sense, family of elliptic operators of positive order
on a compact manifold with boundary and with elliptic boundary conditions
leads to a graph continuous family of regular operators with compact resolvents 
acting between the $L^2$-spaces of the corresponding vector bundles.
See \cite[Appendix A.5]{Pr17} or \cite[Corollary A.6.3]{BBZ}.
A fixed manifold with boundary may be replaced by a family of such manifolds, as above.
Theorems \ref{Thm:RK-H} and \ref{Thm:RKsa-fd} show 
that such a family of unbounded operators may be turned into a Riesz continuous family
by appropriate global trivializations of the corresponding Hilbert bundles,
both in the non-self-adjoint and self-adjoint case.

\sub{Acknowledgments.}
This work was partially supported by ISF grants no. 431/20 and 876/20 
and by the European Research Council (ERC) under the European Union's Horizon 2020 research and innovation programme (grant no.~101001677).
It was done during my postdoctoral fellowship at the Technion -- Israel Institute of Technology; 
the revised version was prepared during my postdoctoral fellowship at the University of Haifa.

I am grateful to N.~V.~Ivanov for his interest in my work and, in particular, for attracting my attention to the work of Dixmier and Douady \cite{DD}. 
The exposition in this paper benefitted a lot from numerous suggestions of the anonymous referee, to whom I am deeply grateful.

\section{Preliminaries}\label{sec:prel}

\uppskip
\sub{Grassmanian.}
The set $\Gr(H)$ of closed subspaces of $H$ is called the \textit{Grassmanian}.
An orthogonal projection may be identified with its range,
which provides a bijection between $\Gr(H)$ and the space of projections $\Proj(H)$.
The Grassmanian is equipped with the topology (and even metric) induced by this bijection
from the norm on $\Proj(H)$.
We will use closed subspaces and projections interchangeably, 
depending on what is more convenient at the moment.

\sub{Canonical decomposition.} 
Recall a standard construction that will be used throughout the paper.
The trivial Hilbert bundle over $\Gr(H)$ with the fiber $H$ is canonically decomposed into the direct sum 
\begin{equation*}
	\Gr(H)\times H = \H^+\oplus\H^- 
\end{equation*}
of two Hilbert bundles, whose fibers over $L\in\Gr$ are 
$\H^+_L = L$ and $\H^-_L = L^\bot$, the orthogonal complement of $L$ in $H$.
This decomposition is locally trivial in the following sense: for every $L_0\in\Gr(H)$ 
there is a continuous map $g\colon B\to\U(H)$, where $B$ is the open ball of radius $1$ around $L_0$, 
such that $L = g_L L_0$ for every $L\in B$.
In terms of projections, this is written as $p = g_p p_0 g_p^*$ for $\norm{p-p_0}<1$,
where $p_0$ is the projection onto $L_0$.
See the proof of \cite[Theorem I.6.32]{Kato} for the construction of such a map $g$.

An important subspace of $\Gr(H)$ is
\[ \Gr\st(H) = \sett{L\in\Gr(H)}{\dim L = \infty = \codim L}. \]
Let $L_0\in\Gr\st(H)$ be a fixed subspace.
The map $\U(H)\to\Gr\st(H)$, $u\mapsto uL_0$ 
is a locally trivial principal bundle over a paracompact base space,
whose structure group $\U(L_0)\times\U(L_0^\bot)$ is contractible.
Therefore, $\U(H)\to\Gr\st(H)$ has a section
\begin{equation}\label{eq:gP*U}
	 g\colon \Gr\st(H)\to\U(H), \txt{so that} L = g_L L_0 \txt{for all} L\in\Gr\st(H).
\end{equation}
It follows that the restrictions of $\H^+$ and $\H^-$ to $\Gr\st$ are trivial Hilbert bundles 
with the structure group $\U$.
In addition, contractibility of the total space $\U(H)$ and the fiber $\U\times\U$ 
implies contractibility of the base space $\Gr\st$. 
See \cite[proof of Lemma 3.6]{ASi}.

\sub{Bounded transform.}
Let us recall the definition and main properties of the bounded transform.
See, e.g., \cite{Schm}, Theorem 7.5 for the detail.

For a regular operator $A\colon H\to H'$, the operator $1+A^*A$ is regular, self-adjoint, and surjective; 
its inverse $(1+A^*A)\inv$ is a positive bounded operator. 
The \emph{bounded transform} (or the Riesz map) 
\[ \f\colon\Reg(H,H')\to\B(H,H'), \quad \f(A) = A(1+A^*A)^{-1/2}, \] 
defines the injective map from the set $\Reg(H,H')$ of regular operators to the closed unit ball 
\[ \D(H,H')=\sett{a\in\B(H,H')}{\norm{a}\leq 1}. \]
The image of $\f$ consists of strict contractions, that is, operators $a\in\D$ such that $1-a^*a$ is injective;
this image is dense in $\D$.
For $a=\f(A)$ the following identities hold: 
\begin{equation}\label{eq:AAaa}
	(1+A^*A)\inv = 1-a^*a \txt{and} A = a(1-a^*a)^{-1/2} .
\end{equation}
In particular, a regular operator $A$ has compact resolvents if and only if $\f(A)$ is essentially unitary.
If a regular operator $A$ is self-adjoint, then so is $\f(A)$; more generally, $\f(A^*) = \f(A)^*$.

\sub{Graphs of operators.}
Let $\Gamma_A\subset H\oplus H'$ denote the graph of a regular operator $A\colon H\to H'$.
There is a canonical isometric isomorphism 
\[ \Gamma_A\to H, \quad \xi\oplus A\xi\mapsto \sqrt{1+A^*A}\,\xi \;\text{ for } \xi\in H. \]
The inverse isomorphism $H\to\Gamma_A$, 
which we consider as an isometry $\w_A\colon H\to H\oplus H'$ with the range $\Gamma_A$,
is given by the formula 
\begin{equation}\label{eq:wA}
 \w_A = (1+A^*A)^{-1/2} \oplus A(1+A^*A)^{-1/2} = (1+A^*A)^{-1/2} \oplus \f(A).  
\end{equation}
Let $\p_A$ denote the orthogonal projection onto $\Gamma_A$,
\begin{equation}\label{eq:pA}
 \p_A = \w_A\w_A^* = \matr{(1+A^*A)^{-1} & (1+A^*A)^{-1}A^* \\ A(1+A^*A)^{-1} & 1-(1+AA^*)^{-1}} \in\Proj(H\oplus H'). 
\end{equation}
In terms of $a=\f(A)$ this is written as follows:
\begin{equation}\label{eq:wa-pa}
  \w_A = \sqrt{1-a^*a} \oplus a, \quad\quad 
	 \p_A = \matr{1-a^*a & \sqrt{1-a^*a}\,a^* \\ a\sqrt{1-a^*a} & aa^*},
\end{equation}
so $\p$ factors through a continuous map $\D(H,H')\to\Proj(H\oplus H')$.
Restricting this map to the subspace $\rR$ of $\D$, we see that the identity map $\rR\to\gR$ is continuous.

\sub{Fredholm operators.}
Similar to the bounded case, a regular operator is called Fredholm 
if its range is closed and its kernel and cokernel are finite-dimensional. 
The bounded transform $\f(A)$ is equal to the composition of the isomorphism $\w_A\colon H\to\Gamma_A$
and the restriction $p'_A\colon \Gamma_A \to H'$ of the orthogonal projection $p'\colon H\oplus H' \to H'$.
The range of $p'_A$ coincides with the range of $A$.
The kernel of $p'_A$ is isomorphic to the kernel of $A$.
Therefore, a regular operator $A$ is Fredholm if and only if $\f(A)$ is Fredholm.

\sub{Graph-to-strong* continuity of the bounded transform.}
By the definition of Riesz topology, the bounded transform $\f\colon\Reg(H,H')\to\D(H,H')$ 
is Riesz-to-norm continuous.
In the proof of Theorem \ref{Thm:Phi} we also need another type of continuity of $\f$.

Let $\sB(H,H')$ denote the space of bounded operators $H\to H'$ equipped with the strong operator topology.
Recall that this topology is defined as the weakest topology for which the map
$\sB(H,H')\to H'$, $a\mapsto a\xi$ is continuous for every $\xi\in H$.
The strong* operator topology on the same set is the weakest topology for which both
$a\mapsto a\xi$ and $a\mapsto a^*\xi$ are continuous for every $\xi\in H$.

\begin{prop}\label{prop:f-g2s}
The bounded transform is graph-to-strong* continuous.
\end{prop}

\proof
Since $\p_{A^*} = 1-J\p_A J^*$, where $J = \smatr{0 & -1 \\ 1 & 0}$, 
the map $A\mapsto A^*$ is graph continuous. 
Since $\f(A)^*=\f(A^*)$, it is enough to prove graph-to-strong continuity of $\f$.

For $A\in\Reg(H,H')$ we denote 
 \begin{equation}\label{eq:aqb}
	 a=\f(A), \quad q = (1+A^*A)^{-1/2} = (1-a^*a)^{1/2}, \quad b = aq. 
 \end{equation}
In these terms, the second equality of \eqref{eq:wa-pa} is written as
\begin{equation}\label{eq:pAa}
  \p_A = \matr{q^2 & qa^* \\ aq & aa^*} = \matr{q^2 & b^* \\ b & aa^*}.
\end{equation}
Since the map $A\mapsto b$ is given by the $(2,1)$-component of $\p_A$, it is graph-to-norm continuous.
Similarly, the map $A\mapsto q$ 
is the square root of the $(1,1)$-component of $\p_A$ and thus is graph-to-norm continuous.

Choose $\xi\in H$ and $A_0\in\Reg(H,H')$.
Let $a_0$, $q_0$, $b_0$ be determined by $A_0$ as in \eqref{eq:aqb}.
We need to show that the map $\gR(H,H')\to H'$, $A\mapsto\f(A)\xi$ is continuous at $A=A_0$.
Let $\eps>0$.
The operator $q_0$ is self-adjoint and injective and thus has dense range, 
so one can find $\eta\in H$ such that $\norm{q_0\eta-\xi}<\eps$.
Let $V$ be a neighborhood of $A_0$ in $\gR(H,H')$ such that 
$\norm{(q-q_0)\eta}<\eps$ and $\norm{(b-b_0)\eta}<\eps$ for every $A\in V$.
Then 
\begin{gather*}
	\norm{a\xi-a_0\xi} \;<\; \norm{\br{a-a_0}q_0 \eta} + 2\eps 
			\;\leq\; \norm{a\br{q_0-q}\eta} + \norm{\br{aq-a_0q_0}\eta} + 2\eps \\
				\leq\; \norm{\br{q-q_0}\eta} + \norm{\br{b-b_0}\eta} + 2\eps < 4\eps
\end{gather*}
for every $A\in V$. 
It follows that $\f$ is graph-to-strong continuous at an arbitrary point $A_0\in\Reg$ and thus on the whole $\Reg$.
\endproof

\begin{rem*}
In fact, this reasoning provides a stronger result. 
It concerns another natural topology on $\Reg(H,H')$ 
induced by the injection $\p\colon\Reg(H,H')\hookto\sP(H\oplus H')$ 
from the strong operator topology on the space of projections. 
On $\Rsa(H)$ it coincides with the strong resolvent topology.
The reasoning in the proof of Proposition \ref{prop:f-g2s} shows that the bounded transform $\f\colon\Reg\to\D$ 
is continuous with respect to this coarser topology on $\Reg$ and the strong* operator topology on $\D$.

For an arbitrary function $f\in C_b(\R)$, a similar reasoning provides continuity of the map $f\colon\Rsa(H)\to\B(H)$ 
with respect to the strong resolvent topology on $\Rsa$ and the strong topology on $\B$.
A different proof of such continuity of $f$ based on the spectral theorem 
is given by \cite{RS}, Theorem VIII.20(b).
Conversely, one can deduce Proposition \ref{prop:f-g2s} from that theorem
taking $f(x)=x(1+x^2)^{-1/2}$ and using the standard trick of passing from $A$ to the self-adjoint operator $\smatr{0 & A^* \\ A & 0}$.
However, the direct proof of the proposition given above and not based on the spectral theorem seems to be more transparent. 
\end{rem*}

\upskip
\sub{Action of the unitary groups.}
The unitary groups $\U(H)$ and $\U(H')$ act on $\Reg(H,H')$ by the right and left multiplication respectively.
These actions preserve the subsets $\RF$ and $\RK$. 

\begin{prop}\label{prop:UR}
The actions of the groups $\U(H)$ and $\U(H')$ (equipped with the norm topology) 
on $\rR(H,H')$ and $\gR(H,H')$ are continuous.
\end{prop}

\proof
For $u\in\U(H')$ and $v\in\U(H)$ we have
\[ \f(uA)=u\f(A), \quad \f(Av) = \f(A)v, \quad \Gamma_{uA} = (1\oplus u)\Gamma_A, \quad 
   \Gamma_{Av} = (v\inv\oplus 1)\Gamma_A,  \]
which proves the proposition.
\endproof

Let $\sU(H)$ denote the unitary group equipped with the strong operator topology. 
It is a topological group, see \cite[Lemma 1.3]{RW}. 
The actions of $\sU(H)$ and $\sU(H')$ on $\rR$ and $\gR$ are not continuous,
but continuity of the action of $\sU$ on $\K$ provides the following result.

\begin{prop}\label{prop:sURK}
The actions of the groups $\sU(H)$ and $\sU(H')$ on $\gR_K(H,H')$ are continuous.
\end{prop}

\proof
The graph of an operator with compact resolvents belongs to the restricted Grassmanian
\begin{equation}\label{eq:PK}
	 \Gr_K(H\oplus H')\cong\PK(H\oplus H') = \sett{p\in\Proj(H\oplus H')}{p-\pinf \text{ is compact}},
\end{equation}
where $\pinf\in\Proj(H\oplus H')$ is the projection onto $0\oplus H'$.
The restriction of the embedding $\p\colon\gR\hookto\Proj$ to operators with compact resolvents
provides the embedding $\p\colon\gR_K\hookto\PK$. 
The left action of $\sU(H')$ on $\gR(H,H')$ lifts to the action on $\Proj(H\oplus H')$:
\begin{equation*}
  (u,p)\mapsto (1\oplus u)p(1\oplus u\inv) = (1\oplus u)(p-\pinf)(1\oplus u)\inv + \pinf.
\end{equation*}
Since the action of $\sU(H\oplus H')$ on $\K(H\oplus H')$ by conjugation is continuous,
the action of $\sU(H')$ on $\PK(H\oplus H')$ is also continuous, as well as its restriction to $\gR_K(H,H')$.
The right action of $\sU(H)$ on $\gR_K$ is continuous by the same argument.
\endproof

\section{Proof of Theorem \ref{Thm:Phi}}\label{sec:twist}\label{sec:Phi}

The idea of the proof can be breifly described as follows.
The graphs $\Gamma_A$ of regular operators $A\colon H\to H'$ form fibers 
of a locally trivial Hilbert bundle $\Gamma$ over $\gR(H,H')$ with the structure group $\U(H)$.
Since $\U(H)$ is contractible, this bundle admits a trivialization $u\colon \gR\times H \isor\Gamma$. 
The composition of the unitary operator $u_A\colon H \isor\Gamma_A\subset H\oplus H'$ with the projection 
$H\oplus H'\to H'$ provides the contraction $\phi_A\colon H\to H'$.
We define $\Phi$ by the formula $\f(\Phi_A) = \phi_A$
and show that this map satisfies all conclusions of the theorem.

\sub{Isometries.}
Let 
\[ \I = \sett{u\in\B(H, H\oplus H')}{u^*u = \Id_{H}} \] 
be the space of isometries (i.e., isometric embeddings) $H\to H\oplus H'$ equipped with the usual norm topology.
The canonical map
\begin{equation*}
	\ppi\colon\I\to\Gr(H\oplus H'), \quad u\mapsto \Ran(u). 
\end{equation*}
takes an isometry $u$ to its range; the corresponding projection is $uu^*$.
The image $\ppi(\I)$ of this map consists of infinite-dimensional subspaces of $H\oplus H'$.
The group $\U(H)$ acts on $\I$ by the right multiplication.
Kuiper's theorem \cite{Kui} together with the same standard arguments as were used 
in the beginning of Section \ref{sec:prel} show that 
$\ppi\colon\I\to\ppi(\I)$ is a principal $\U(H)$-bundle admitting a continuous section (and thus is trivial).

\sub{Proof of Theorem \ref{Thm:Phi}.}
The graphs of regular operators have infinite dimension, so $\gR\subset\ppi(\I)\subset\Gr(H\oplus H')$.
For the rest of the proof we fix a section of $\ppi\colon\I\to\ppi(\I)$ 
and denote by $u\colon\gR\to\I$ its restriction to $\gR$.
Then $u$ is a graph-to-norm continuous map taking a regular operator $A$ 
to an isometry $u_A\colon H\to H\oplus H'$ with the range $\Gamma_A$.

Let $p'\colon H\oplus H' \to H'$ be the projection onto the second summand.
The composition
\[ \phi_A = p'\cdot u_A\colon H\to H' \]
is a contraction, and the corresponding map $\phi\colon\gR(H,H')\to\D(H,H')$ is continuous.

The maps defined above are not canonical: they depend on the choice of a section $u$.
On the other hand, formula \eqref{eq:wA} provides a canonical isometry $\w_A\colon H\to H\oplus H'$ 
with the range $\Gamma_A$:
\[ \w_A = (1+A^*A)^{-1/2} \oplus \f(A) = \sqrt{1-a^*a} \oplus a, \txt{where} a=\f(A). \]
The corresponding contraction $p'\cdot \w_A\colon H\to H'$ coincides with the bounded transform $\f(A)$ of $A$.
In contrast with $\phi$, it is \textit{not} continuous as a map $\gR\to\D$.

Two isometries $u_A, \w_A\colon H\to H\oplus H'$ have the same range $\Gamma_A$ 
and thus differ by a unitary operator 
\begin{equation*}
	v_A = \w_A^*\cdot u_A\in\U(H), \txt{so that} u_A = \w_A\cdot v_A.
\end{equation*}
We claim that the map $v\colon\Reg(H,H')\to\U(H)$ defined by this formula 
satisfies conclusions of the theorem.
Indeed, for $\Phi_A=A\cdot v_A$ we have
\[ \f(\Phi_A) = \f(A\cdot v_A) = \f(A)\cdot v_A = p'\cdot \w_A\cdot v_A = p'\cdot u_A = \phi_A, \]
so $\f\circ\Phi = \phi\colon\Reg(H,H')\to\B(H,H')$ is graph-to-norm continuous 
and thus $\Phi$ is graph-to-Riesz continuous.
By the construction, $u\colon\Reg\to\I$ is graph-to-norm (and thus also Riesz-to-norm) continuous.
The coisometry $\w_A^*\colon H\oplus H'\to H$ is given by the formula 
\[ \w_A^*(\xi\oplus\eta) \;=\; (1+A^*A)^{-1/2}\;\xi + \f(A)^*\,\eta \;=\; \sqrt{1-a^*a}\;\xi + a^*\,\eta  
   \;\txt{for} \xi\in H, \; \eta\in H'. \]
Clearly, the maps $A\mapsto \w_A^*$ and $A\mapsto v_A = \w_A^*\cdot u_A$ are Riesz-to-norm continuous.

The first component of $\w^*$, $A\mapsto (1+A^*A)^{-1/2}$, 
is the square root of the $(1,1)$-component of $\p_A$ and thus is even graph-to-norm continuous.
The second component $A\mapsto \f(A)^*$ 
is graph-to-strong continuous by Proposition \ref{prop:f-g2s}.
Therefore, the map $A\mapsto \w_A^*$ is graph-to-strong continuous.
The multiplication map $(b,c)\mapsto b\cdot c$ is strongly continuous on uniformly bounded subsets of operators,
and the norms of $\w_A^*$ and $u_A$ are bounded by $1$.
It follows that the map $A\mapsto v_A = \w_A^*\cdot u_A$ is graph-to-strong continuous.
This completes the proof of the theorem.
\endproof

\section{Applications of Theorem \ref{Thm:Phi}}\label{sec:appl1}

\upskip
\sub{Proof of Theorem \ref{Thm:RKF}.}
In order to prove the first part of the theorem, we need to show that $\Phi$ is homotopic to the identity 
as a self-map of both $\rR_K$ and $\gR_K$.
For the second part we need to prove a similar result for $\rR_F$ and $\gR_F$.

By Kuiper's theorem, the group $\U(H)$ is contractible in the norm topology, 
so the map $v\colon\U\to\U$ is homotopic to the identity map $\U\to\U$.
Since the action of $\U(H)$ on $\rR(H,H')$ is continuous,
a homotopy between these two maps provides a homotopy between 
$\Phi\colon\rR_K\to\rR_K$ and $\Id\colon\rR_K\to\rR_K$, 
and a similar homotopy for $\rR_F\to\rR_F$.

By a result of Dixmier and Douady \cite[Lemma 3]{DD}, the group $\sU(H)$ is contractible,
so the map $v\colon\sU\to\sU$ is homotopic to the identity map $\sU\to\sU$.
The action of $\sU(H)$ on $\gR_K(H,H')$ 
is continuous (see Proposition \ref{prop:sURK}),
so a homotopy between these two maps provides a homotopy 
between the maps $\Phi\colon\gR_K\to\gR_K$ and $\Id\colon\gR_K\to\gR_K$.
Therefore, the restriction of $\Phi$ to the subspace of operators with compact resolvents 
is homotopy inverse to the identity map $\rR_K\to\gR_K$.

In order to finish the proof of the Fredholm part, 
we need to show that the maps $\Phi\colon\gR_F\to\gR_F$ and $\Id\colon\gR_F\to\gR_F$ are homotopic.
This part is more delicate because the action of $\sU$ on $\gR_F$ is not continuous,
so we cannot directly apply a contraction of $\sU$ for this purpose.
However, we can use a part of \cite[Theorem A]{Pr21} 
which states that the embedding $\gR_K\hookto\gR_F$ is a homotopy equivalence.
We already proved that the maps $\Phi\colon\gR_K\to\gR_K$ and $\Id\colon\gR_K\to\gR_K$ are homotopic; 
it follows that $\Phi\colon\gR_F\to\gR_F$ and $\Id\colon\gR_F\to\gR_F$ are also homotopic.
Therefore, the restriction of $\Phi$ to the subspace of Fredholm operators 
is homotopy inverse to the identity map $\rR_F\to\gR_F$.
\endproof

\sub{Proof of Theorem \ref{Thm:RK-H}.}
Suppose that a global trivialization of $\H'$ is fixed.
Choose an arbitrary trivialization $\tau\colon\H\to X\times H$ of $\H$.
In these trivializations $\A$ becomes a graph continuous map $X\to\RK(H,H')$.
Let $v\colon\gR(H,H')\to\sU(H)$ be a map provided by Theorem \ref{Thm:Phi}.
Then the map $X\to\sU(H)$, $x\mapsto v(\A_x)$ is continuous, so 
\[ \tau_{v}\colon(x,\xi)\mapsto(x, v(\A_x)\xi) \] 
is an automorphism of the trivial Hilbert bundle $X\times H$.
The inverse automorphism $\tau_{v}\inv$ takes $\A$ to a Riesz continuous map 
\[ \A'\colon X\to\RK(H,H'), \quad \A'_x = \A_x\cdot v(\A_x). \]
It follows that the composition $\tau_{v}\inv\circ\tau$ satisfies conclusion of the theorem.

If instead a global trivialization of $\H$ is fixed, 
then one can pass to the adjoint operators, swapping $\H$ and $\H'$.
Since the map $A\mapsto A^*$ is both graph and Riesz continuous, 
it follows from the first part of the proof 
that there is a global trivialization of $\H'$ taking $\A$ to a Riesz continuous map.
\endproof

\sub{$K^0$ index for graph continuous families.}
By a classical result of Atiyah and J\"anich, 
the space $\BF(H)$ of bounded Fredholm operators is a classifying space for the functor $K^0$.
See \cite{Atiyah}, Theorem A1.
For an arbitrary parameter space $X$,
$K^0(X)$ might be \textit{defined} as the set $[X,\BF]$ of homotopy classes of maps $X\to\BF(H)$,
with the group structure on $[X,\BF]$ introduced in the standard way.
The class of a continuous map $\alpha\colon X\to\BF$ in $K^0(X)$ is called the \textit{index} of $\alpha$:
\[ \ind(\alpha) = [\alpha]\in[X,\BF]=K^0(X). \]
It can be easily seen, using linear homotopies, that the bounded transform $\f\colon\rR_F\to\BF$,
from regular Fredholm operators to bounded Fredholm operators,
is homotopy inverse to the embedding $\BF\hookto\rR_F$,
so both these maps are homotopy equivalencies.
Therefore, the index of a Riesz continuous map $\A\colon X\to\rR_F$ 
is naturally defined as the index of its bounded transform: 
\begin{equation}\label{eq:ind-rR}
	\ind(\A) = \ind(\f\circ\A)\in[X,\BF]=K^0(X) \txt{for} \A\colon X\to\rR_F. 
\end{equation}

Let now $\A\colon X\to\gR_F$ be a \emph{graph} continuous map.
The space $\gR_F$ is also a classifying space for the functor $K^0$,
see \cite{Jo}, Theorem 3.5(i) or \cite{Pr21}, Theorem A.
Hence the homotopy class of $\A$ determines an element of $K^0(X)$ which can be called the index of $\A$:
\begin{equation}\label{eq:ind-gR}
	\ind(\A) = [\A]\in[X,\gR_F]=K^0(X). 
\end{equation}
Theorems \ref{Thm:Phi} and \ref{Thm:RKF} of the present paper allow to relate the index of $\A$ 
given by \eqref{eq:ind-gR} to the Fredholm index of bounded operators as follows.
The composition of $\A$ with a map $\Phi\colon\gR\to\rR$ provided by Theorem \ref{Thm:Phi} is Riesz continuous,  
so $\f\circ\Phi\circ\A$ is a norm continuous map $X\to\BF$.
Theorem \ref{Thm:RKF} implies
\begin{equation}\label{eq:ind-gr}
	\ind(\A) = \ind(\Phi\circ\A) = \ind(\f\circ\Phi\circ\A) \in [X,\BF] = K^0(X) \txt{for} \A\colon X\to\gR_F, 
\end{equation}
which can be used as the definition of the index of $\A$.
It follows from Theorem \ref{Thm:RKF} that the homotopy class of $\Phi\colon\gR_F\to\rR_F$ 
is independent of the choice of $v$, so $\ind(\f\circ\Phi\circ\A)$ is also independent of this choice.

Similarly, Theorem \ref{Thm:RKF} implies that
for a Riesz continuous map $\A\colon X\to\RF$ 
(in particular, for a norm continuous map $X\to\BF$)
the two definitions of the index given by \eqref{eq:ind-rR} and \eqref{eq:ind-gR} coincide.

Recently N. Ivanov \cite{I1, I2} developed a new approach to the family index.
He introduced the notion of a \q{Fredholm family} and defined a $K^0(X)$-valued index, 
which we will denote by $\Ind$, 
for Fredholm families parametrized by points of a paracompact space $X$.
See \cite{I2}, Section 3.
A Fredholm family is a family of Fredholm operators $\A_x$
such that the finite-dimensional subspaces $\one_{[0,\eps]}(|\A_x|)$ and 
$\one_{[0,\eps]}(|\A^*_x|)$, as well as the restriction of $\A_x$ to $\one_{[0,\eps]}(|\A_x|)$,
depend continuously on $x$ for an appropriate locally chosen $\eps>0$,
without any additional continuity assumptions.
Ivanov's index $\Ind$ is defined both for operators acting on a fixed Hilbert space and on fibers of a Hilbert bundle;
it is invariant under isomorphisms of Hilbert bundles and under passing to the bounded transform.
By the construction, for a norm continuous map $\A\colon X\to\BF$ Ivanov's index coincides 
with the usual index $\ind(\A)\in[X,\BF] = K^0(X)$.

Let now $\A\colon X\to\gR_F$ be a graph continuous map from a paracompact space $X$.
Then $\A$ is a Fredholm family in Ivanov's sense, 
and the properties of $\Ind$ mentioned above imply the equality 
\[ \Ind(\A) = \Ind(\f\circ\Phi\circ\A) = \ind(\f\circ\Phi\circ\A) . \]
Therefore, our definition of the index of $\A$ provides the same element of $K^0(X)$ as Ivanov's definition.

\begin{rem*}
Our proof of Theorem \ref{Thm:Phi} shows that the definition of the index of $\A\colon X\to\gR_F$
given by the formula $\ind(\A) = \ind(\f\circ\Phi\circ\A)$
can be equivalently stated in terms of the graphs $\Gamma_x$ of operators $\A_x$,
as the index of the family of projections $\Gamma_x\to H'$.
For this, one needs to trivialize the Hilbert bundle $\Gamma$, formed by subspaces $\Gamma_x\subset H\oplus H'$,
in the category of bundles with the structure group $\U$.
Note that such a trivialization is not canonical; it is defined only up to a continuous map $X\to\U(H)$.
However, change of a trivialization does not affect the homotopy type of the corresponding family
of bounded Fredholm operators, since the unitary group $\U(H)$ is contractible.
\end{rem*}

\section{Polarizations and proof of Theorem \ref{Thm:Gr}}\label{sec:polar}

\upskip
\sub{Polarizations of a Hilbert space.}
A decomposition of a Hilbert space $H$ into the orthogonal sum $H = H^+\oplus H^-$ 
of two closed infinite-dimensional subspaces is called a \textit{polarization} of $H$.
It is determined by the subspace $H^+\in\Gr\st(H)$ 
or equivalently by the orthogonal projection onto $H^+$.
We denote by $\Po(H)\cong\Gr\st(H)$ the subspace of $\Proj(H)$ consisting of projections of infinite rank and corank
and will identify polarizations with the corresponding subspaces/projections.

There are two natural equivalence relations for polarizations of $H$, which we call \q{finite} and \q{compact}.
Two polarizations $p$ and $q$ of $H$ are said to be compactly compatible if $p-q$ is a compact operator;
they are said to be finitely compatible if $p-q$ is a finite rank operator
(that is, the corresponding subspaces are commensurable).

The compact restricted Grassmanian corresponding to a polarization $p$ 
is defined as the space of all polarizations compatible with $p$: 
\[ \Gr_K(H;p) = \sett{q\in (H)}{q-p\in\K(H)}. \]
It is equipped with the norm topology induced by the natural inclusion $\Gr_K(H;p)\hookto\Gr(H)\cong\Proj(H)$, 
or equivalently by the inclusion $\Gr_K(H;p)\hookto\K(H)$, $q\mapsto q-p$.
The finite restricted Grassmanian $\Gr_f(H;p)$ is defined in a similar manner; it is a subspace of $\Gr_K(H;p)$.
The restricted Grassmanian depends only on the corresponding equivalence class of the projection $p$.

\sub{Subbundles of a Hilbert bundle.}
We will consider only locally trivial Hilbert bundles with a separable fiber $H$. 
However, we want to deal with their subbundles which are not necessarily locally trivial.
Such subbundles arise, for example, as spectral sections of graph continuous families of Fredholm operators.

In general, a subbundle $\E$ of $\H$ is given by a family of closed subspaces $\E_x\subset\H_x$ 
which is continuous in an appropriate sense.
We will use the following notion of continuity:
we say that $\E$ is a \emph{subbundle} of $\H$ if $\E$ is the range of a bundle projection, 
that is, a continuous bundle map $p\colon\H\to\H$ such that $p_x\in\Proj(\H_x)$ for every $\xX$.
In this case $\H\ominus\E = \Ran(1-p)$ is also a subbundle. 

Equivalently, for every $\xX$, some (and then every) local trivialization of $\H$ 
over a neighborhood $Y$ of $x$ takes $p$ to a continuous map $Y\to\sP(H)$, 
where $\sP(H)$ denotes the space of projections equipped with the strong operator topology.
In particular, subbundles of the trivial Hilbert bundle $H_X=X\times H$ are in one-to-one correspondence 
with continuous maps $X\to\sP(H)$. 

A family $\E$ of closed subspaces $\E_x\subset\H_x$ is a subbundle of $\H$ in our sense 
if and only if both $(\E_x)$ and $(\H_x\ominus\E_x)$ are continuous subfields of the field $(\H_x)$ of Hilbert spaces in the sense of Dixmier and Douady \cite{DD}. 
See \cite{DD}, Proposition 16.

We say that two subbundles of $\H$ given by bundle projections $p$ and $q$
are (compactly) compatible if $p-q$ is a norm continuous family of compact operators;
they are finitely compatible if additionally all operators $p_x-q_x$ have finite rank.

\sub{Polarizations of a Hilbert bundle.}
A natural forgetful functor from finite to compact polarizations preserves the notion of an adapted trivialization.
Therefore, the \q{finite} analogues of Theorems \ref{Thm:Gr} and \ref{Thm:Gr-fd}  
follow immediately from these theorems. 
Because of that, we will deal only with compact polarizations in the rest of the paper 
and will mostly omit the adjective \q{compact} before \q{polarization}.
The changes required to obtain the finite version of our statements are obvious.

Locally, a polarization of $\H$ is determined by a decomposition of $\H$, over some open subset $Y\subset X$, 
into the orthogonal sum $\H = \H^+\oplus \H^-$ of two subbundles of infinite rank.
Equivalently, it is determined by the subbundle $\H^+$ of $\restr{\H}{Y}$ 
or by the bundle projection with the range $\H^+$. 
By the definition, two compatible subbundles of $\H$ over $Y$ determine the same local polarization.

Globally, a polarization of $\H$ is determined by an open covering $(X_i)$ of $X$ 
and subbundles $\H^+_i$ of $\restr{\H}{X_i}$ as above such that 
$\H^+_i$ and $\H^+_j$ are compatible over $X_i\cap X_j$ for every $i,j$.
Following Ivanov \cite{I2}, we call such a collection $(X_i,\H^+_i)$ a \textit{polarization atlas}.
Two polarization atlases are said to be equivalent if their union is also a polarization atlas. 
A polarization of $\H$ is an equivalence class of polarization atlases.

A polarization $\Pi$ of $\H$ determines the restricted Grassmanian bundle $\Gr_K(\H;\Pi)$ in an obvious way.
This is one of the motivations behind the notion of a polarization.

Two polarizations of $\H$ are said to be homotopic if they are restrictions 
to $X\times\set{0}$ and $X\times\set{1}$ of some polarization of the Hilbert bundle 
$\H\times[0,1]\to X\times[0,1]$.

\sub{Adapted trivializations.}
A trivialization $\H\to X\times H$ is said to be \textit{adapted} to a subbundle $\E$ of $\H$ 
if it takes the family of subspaces $\E_x\subset\H_x$ to a continuous map $X\to\Gr(H)$,
or equivalently takes the corresponding family of projections to a norm continuous map $X\to\Proj(H)$.
If a trivialization of $\H$ is adapted to $\E$, then it is adapted to every subbundle compatible with $\E$,
since the difference $p-q$ of corresponding bundle projections is norm continuous in every trivialization.

Let $\Pi$ be a polarization given by an atlas $(X_i,\E_i)$.
We say that a trivialization of $\H$ is \textit{adapted} to $\Pi$ 
if its restriction to $X_i$ is adapted to $\E_i$ for every $i$. 
The notion of an adapted local trivialization is defined in an obvious way.

We call a subbundle $\E$, resp. a polarization $\Pi$ of $\H$ \textit{tame} 
if it admits local adapted trivializations 
(that is, for every point $x\in X$ there is a local trivialization of $\H$ adapted to $\E$, resp. $\Pi$ 
over some neighborhood of $x$).

\sub{Trivializations of $H_X$ adapted to $E_X$.}
As a first step to the proof of Theorem \ref{Thm:Gr},
we describe trivializations of the trivial Hilbert bundle $H_X=X\times H$ 
adapted to its trivial subbundle $E_X=X\times E$,
where $E$ is a fixed closed subspace of $H$.

Trivializations of $H_X$ are given by continuous maps $f\colon X\to\sU(H)$.
Such a trivialization is adapted to $E_X$ if and only if the map $X\to\Gr(H)$ 
taking $x$ to $f(x)E$ is continuous.
Therefore, trivializations of $H_X$ adapted to $E_X$ are in one-to-one correspondence 
with continuous maps $X\to T_E$, where the space $T_E$ is obtained from  $\sU(H)$ 
by refining its topology so that the map
\begin{equation}\label{eq:EsUP}
  \iota_E\colon T_E\to \Gr(H), \quad u\mapsto uE
\end{equation}
becomes continuous.
In other words, $T_E$ has the topology induced by the injection 
\begin{equation*}
  T_E\hookto \sU(H)\times\Gr(H), \quad u\mapsto(u, uE). 
\end{equation*}
Two maps $f,g\colon X\to T_E$ are homotopic if and only if 
the corresponding trivializations are homotopic in the class of trivializations adapted to $E_X$.

\begin{lem}\label{lem:TE}
  If $E\in\Gr\st(H)$, then $T_E$ is contractible.
\end{lem}

\proof
For $E\in\Gr\st$ the map \eqref{eq:EsUP} takes $T_E$ to $\Gr\st$.
We will show that $\iota_E\colon T_E\to \Gr\st$ is a trivial bundle and construct its trivialization.

The base space $\Gr\st$ has a distinguished point $E\in\Gr\st$.
The fiber of $\iota_E$ over this point is the subspace of $\sU(H)$ 
consisting of unitary operators $u\colon H\to H$ preserving $E\subset H$.
Denote this subspace by $\sU_E$.
Let $g\colon\Gr\st\to\U(H)$ be a continuous map such that 
$L = g_L E$ for all $L\in\Gr\st$, as in \eqref{eq:gP*U}.
Consider the map
\begin{equation*}
  \phi\colon \sU_E\times\Gr\st \;\to\; T_E, \quad \phi(v,L) = g_L\cdot v.
\end{equation*}
The composition $\iota_E\circ\phi\colon\sU_E\times\Gr\st \to \Gr\st$ 
is the projection of the product on the second factor, as can be seen from
\begin{equation}\label{eq:v2UE}
	(g_L v)E = g_L E = L \txt{for every} v\in \sU_E.
\end{equation}
In particular, $\iota_E\circ\phi$ is continuous.
The composition of $\phi$ with the identity map $T_E\to\sU$ is given by $(v,L)\mapsto g_Lv$
and is also continuous.
Therefore, $\phi$ itself is continuous.

Equality \eqref{eq:v2UE} allows to recover $L$ and $v$ from $u=g_Lv$ by taking $L=uE$ and $v=g_{L}\inv u$.
Since $vE = g_{L}\inv uE = g_L\inv L = E$, the unitary operator $v=g_{L}\inv u$ belongs to $\sU_E$,
and we obtain the continuous map 
\begin{equation*}
 \psi\colon T_E\to \sU_E\times\Gr\st, \quad \psi(u) = \br{g_{uE}\inv\cdot u,\, uE}. 
\end{equation*}
By the construction, $\psi\circ\phi(v,L)=(v,L)$ for every $(v,L)\in\sU_E\times\Gr\st$.
Conversely,
\[ \phi\circ\psi(u) = \phi\br{g_{uE}\inv\cdot u,\, uE} = g_{uE}\cdot g_{uE}\inv\cdot u = u \txt{for every} u\in T_E. \]
Therefore, $\psi$ is inverse to $\phi$ and thus $\phi$ is a homeomorphism.

The domain of $\phi$ is the product of two factors. 
The second factor $\Gr\st$ is contractible, see \cite[proof of Lemma 3.6]{ASi}.
The first factor $\sU_E = \sU(E)\times\sU(E^\bot)$ is the product of two copies of $\sU$.
By \cite[Lemma 3]{DD}, $\sU$ is contractible, so $\sU_E$ is also contractible.
Therefore, the product $\sU_E\times\Gr\st\cong T_E$ is contractible.
\endproof

\sub{Trivializations of $\H$.}
Let $\H$ be a locally trivial Hilbert bundle and 
$H_X$ be the trivial Hilbert bundle with the fiber $H$ over a topological space $X$. 
Let $\U(\H,H_X)$ be the fiber bundle associated with $\H$ 
whose fiber over $x$ is $\sU(\H_x,H)$.
This is a locally trivial principal bundle over $X$ with the structure group $\sU(H)$
acting on fibers by the left multiplication.
Sections of $\U(\H,H_X)$ are Hilbert bundle isomorphisms $\H\to H_X$, that is, trivializations of $\H$.
A trivialization $\H\to H_X$ takes $\U(\H,H_X)$ to the trivial bundle $\U(H_X,H_X)=X\times\sU(H)$.

\sub{Trivializations adapted to $\E$.}
Let $\E$ be a subbundle of $\H$. 
Consider the space $\T=\T_\E$ which coincides with $\U(\H,H_X)$ as a set, 
but has finer topology determined by $\E$.
Namely, we equip $\T$ with the weakest topology such that both the identity map 
$\iota\colon\T \to \U(\H,H_X)$ and the map
\begin{equation}\label{eq:T2Gr}
	\iota_\E\colon\T \to \Gr(H) \cong \Proj(H), 
	\quad \iota_\E(x,u)\ = u(\E_x) \txt{for} \xX \txt{and} u\in\U(\H_x,H)
\end{equation}
are continuous. 
Then the composition $\T\to\U(\H,H_X)\to X$ is also continuous, and we consider $\T$ as a space over $X$.
By the construction, sections of $\T\to X$ are in one-to-one correspondence 
with trivializations of $\H$ adapted to $\E$.

\begin{lem}\label{lem:EF-comp}
  If $\E$ and $\F$ are compatible subbundles of $\H$, 
	then $\iota_\F$ and $\iota_\E$ induce the same refinement of the original topology on $\U(\H,H_X)$.
	In other words, $\T_\E$ depends only on the equivalence class of compatible subbundles $\E$.
\end{lem}

\proof
Let $p_x$ and $q_x$ be projections onto $\E_x$ and $\F_x$ respectively.
By the definition of compatible subbundles, 
the difference $k_x = q_x-p_x$ is a norm continuous family of compact operators,
so $\iota_\F$ and $\iota_\E$ considered as maps from $\U(\H,H_X)$ to $\Proj(H)$ differ by a continuous map 
\[ \U(\H,H_X)\to\K(H), \quad (x,u)\mapsto \iota_\F(x,u)-\iota_\E(x,u) = uq_xu^* - up_xu^* = uk_xu^*. \] 
Hence $\iota_\F$ and $\iota_\E$ induce the same refinement of the original topology on $\U(\H,H_X)$.
\endproof

\begin{lem}\label{lem:T-loc}
  If $\E$ is tame, then $\T_\E\to X$ is a locally trivial bundle, 
	whose fiber $T=T_E$ is the space defined before Lemma \ref{lem:TE}.
\end{lem}

\proof
Let $x_0\in X$.
Since $\E$ is tame, there is a trivialization $\tau'$ of $\H$ adapted to $\E$ over a neighborhood $X'$ of $x_0$.
This trivialization takes $\E$ to the subbundle $\E'$ given by a \emph{continuous} map $X'\to\Gr(H)$.
Replacing $X'$ by a smaller neighborhood if necessary, 
we can find a map $g\colon X'\to\U(H)$ such that $g_x\E'_x=E$, where $E=\E'_{x_0}$.
The composition of $\tau'$ with the automorphism of the trivial Hilbert bundle $H_{X'}$ provided by $g$
determines the new trivialization $\tau$ which takes $\E_x$ to the subspace $E$ of $H$ independent of $x$
and thus takes $\E$ to the constant subbundle $E_{X'}$ of $H_{X'}$.

Such a trivialization $\tau$ provides a bundle isomorphism $\restr{\T_\E}{X'}\to\T'$,
where the bundle $\T'$ is build (in the way described above) 
from the trivial bundle $H_{X'}$ and its subbundle $E_{X'}$.
By the definition, the total space $\T'$ is obtained from the product $X'\times \sU(H)$
by refining its topology so that the map	$\T'\to \Gr(H)$, $(x,u)\mapsto uE$ becomes continuous.
But this map does not actually depend on $x$, so $\T'$ coincides with the trivial bundle $X'\times T_E\to X'$.
\endproof

\sub{Trivializations adapted to $\Pi$.}
Let $\Pi$ be a polarization of $\H$ given by an atlas $(X_i,\E_i)$, $i\in I$.
Consider the space $\T_\Pi$ which coincides with $\U(\H,H_X)$ as a set 
and is equipped with the weakest topology such that the identity map 
\begin{equation}\label{eq:TPi-TEi}
	\restr{\T_\Pi}{X_i}\to\T_{\E_i} 
\end{equation}
is continuous for every $i\in I$.
By Lemma \ref{lem:EF-comp} the restrictions of $\T_{\E_i}$ and $\T_{\E_j}$ to $X_i\cap X_j$ coincide, 
so \eqref{eq:TPi-TEi} is actually a homeomorphism.
Moreover, this definition does not depend on the choice of an atlas of $\Pi$.

The composition $\T_\Pi\to\U(\H,H_X)\to X$ is continuous, and we consider $\T_\Pi$ as a space over $X$.
By the construction, (local or global) sections of $\T_\Pi\to X$ are exactly 
(local or global) trivializations of $\H$ adapted to $\Pi$.
In particular, $\T_\Pi\to X$ admits local sections if and only if $\Pi$ is tame.

\sub{Proof of Theorem \ref{Thm:Gr}.}
Let $\Pi$ be a tame polarization of a locally trivial Hilbert bundle $\H$ over $X$.
Then, by Lemma \ref{lem:T-loc}, $\T_\Pi\to X$ is a locally trivial bundle with the fiber $T_E$, $E\in\Gr\st(H)$.
This fiber is contractible by Lemma \ref{lem:TE}.
If additionally $X$ is paracompact, then
$\T_\Pi\to X$ admits global sections, every two sections are homotopic, 
and a section over a closed subspace $Y\subset X$ can be extended to a global section \cite[Lemma 4]{DD}.
Since sections of $\T_\Pi$ correspond to adapted trivializations, this proves the theorem.
\endproof

\section{Proof of Theorem \ref{Thm:Gr-fd}}\label{seq:fcd}

Let $\dim Z$ denote the covering dimension of a paracompact space $Z$.
Recall that a closed subset $Y$ of a paracompact space $Z$ 
is itself paracompact \cite[Corollary 5.1.29]{Eng}, normal \cite[Theorem 5.1.5]{Eng},
and $\dim Y \leq \dim Z$ \cite[Theorem 7.1.8]{Eng}.

\sub{Finite covering dimension.}
Our proof of Theorem \ref{Thm:Gr-fd} is based on Theorem \ref{thm:Thm5-DD} below.
Theorem \ref{thm:Thm5-DD} follows from a theorem of Dixmier and Douady \cite[Theorem 5]{DD}.
For the convenience of the reader, we include here a proof of Theorem \ref{thm:Thm5-DD}; 
it is based on the same ideas as the proof of \cite[Theorem 5]{DD}, 
but is stated in terms of operators in a fixed Hilbert space instead of 
\q{continuous fields of Hilbert spaces} of Dixmier--Douady.

Let $E$ be a subspace of $H$ of infinite dimension and codimension and 
$p_E\in\Po(H)$ be the orthogonal projection onto $E$.
In terms of projections, the map $\U(H)\to\Gr\st(H)$, $u\mapsto uE$ 
is written as $\U(H)\to\Po(H)$, $u\mapsto up_Eu\inv$.
This map has a global section, as explained in the beginning of Section \ref{sec:prel}.
In contrast with this, the continuous map 
\begin{equation*}
	\sU(H)\to\sP\st(H) 
\end{equation*}
given by the same formula $u\mapsto up_Eu\inv$ has no sections, even locally.
However, the following result holds.

\begin{thm}\label{thm:Thm5-DD}
	Let $X$ be a paracompact space of finite covering dimension.
	Then every continuous map $p\colon X\to\sP\st(H)$ can be lifted to a map $u\colon X\to\sU(H)$.
	In other words, every subbundle of $H_X$ with infinite rank and corank is trivial.
\end{thm}

\proof
Let $S(H)$ be the unit sphere in $H$. 
For $q\in\Proj(H)$ and $F=\Ran q$ we denote the unit sphere in $F$ by $S(q)=S(F)$.

The map $p\colon X\to\sP\st(H)$ determines a bundle projection $H_X\to H_X$ of the trivial bundle $H_X$, 
which we also denote by $p$.
Let $\E$ be the range of $p$ (it is a subbundle of $H_X$) and 
\[ S(p) = S(\E) := (X\times S(H))\cap\E = \sett{(x,e)}{x\in X,\; e\in S(\E_x) } \]
be the bundle of unit spheres in $\E$.

\begin{lem}\label{lem:DD}
	Let $\eta$ be a continuous section of $\E\to X$ and $\eps>0$.
	Then there is a continuous section $\xi$ of $S(\E)\to X$ such that 
	the distance from $\eta(x)$ to the linear span $\bra{\xi(x)}$ of $\xi(x)$ is less than $\eps$ for every $x\in X$.
\end{lem}

\proof
For a non-zero vector $v\in H$, we will denote by $\s(v) = \norm{v}\inv\cdot v \in S(H)$
the corresponding unit vector.

Consider the family $\S = \sett{S(q)}{q\in\Po(H)}$ of closed subsets of $H$. 
Let us check that the composition $X\to\sP\st\to\S$ 
taking $x\in X$ to $S(p_x)\in\S$ satisfies all four assumptions 
of Michael's Selection Theorem \cite[Theorem 1.2]{Mich}.

(1) $X$ is paracompact and has finite covering dimension.

(2) Every element of the family $\S$ is contractible. 
Indeed, the unit spheres $S(q)$ are contractible by Kakutani's theorem \cite[Theorem 3]{Kaku}.

(3) Let us show that the family $\S$ is equi-$LC^n$ in terms of \cite{Mich} for every $n$.
In fact, the stronger condition holds: 
for every  open ball $V$ of radius $<1$ with the center in $v_0\in S(H)$ and for every $q\in\Po$, 
the intersection $V\cap S(q)$ is either empty or contractible. 
Indeed, if $V\cap S(q)$ is non-empty, then it coincides with $V'\cap S(q)$, 
where $V'$ is an open ball of the same or smaller radius with the center in $v'_0 = \s(qv_0)\in S(q)$. 
The formula $(v,t)\mapsto \s((1-t)v+tv'_0)$ provides a contraction of $V'\cap S(q)$. 
Therefore, $\S$ is equi-$LC^n$.

(4) Let us show that the map $\sP\st(H)\to\S$, $q\mapsto S(q)$ is lower semi-continuous,
that is, for every open $V\subset H$ the set $U=\sett{q\in\Po}{S(q)\cap V\neq\emptyset}$ is open in $\sP\st$.
Let $q_0\in U$ and $v_0\in S(q_0)\cap V$.
Then $W = \sett{q\in\Po}{\norm{(q-q_0)v_0}<1}$ is a neighborhood of $q_0$ in $\sP\st$ 
and the map $f\colon W\to H$, $f(q)=\s(qv_0)\in S(q)$ is continuous and 
has the prescribed value $v_0$ at $q_0$.
Thus the inverse image $f\inv(V)\subset W$ is open in $\sP\st$, contains $q_0$, and is contained in $U$.
Therefore, $U$ is open in $\sP\st$.

We see that the map $X\to\S$, $x\mapsto S(p_x)\in\S$ 
satisfies assumptions of \cite[Theorem 1.2]{Mich}
and thus every section of $S(\E)$ over a closed subset $Y\subset X$ 
can be extended to a global section of $S(\E)$.

Let $Y=\sett{x\in X}{\norm{\eta(x)}\geq\eps}$.
Then the vector field $\xi\colon Y\to S(H)$, $\xi(x) = \s(\eta(x))$ 
is a section of $S(\E)$ over $Y$.
Let $\xi$ be an extension of this section to the whole $X$.
By the construction, the distance from $\eta(x)$ to $\bra{\xi(x)}$ is zero for $x\in Y$
and less than $\eps$ for $x\notin Y$, so $\xi$ is a required section.
\endproof

\sub{Continuation of the proof of Theorem \ref{thm:Thm5-DD}.}
We aim to construct continuous maps $\xi_n\colon X\to S(H)$, $n\in\N$ such that 
$\xi_1(x),\xi_2(x),\ldots$ is an orthonormal basis of $\E_x$ for every $x\in X$.
To this end, fix a countable dense subset $\bar{\eta}_1, \bar{\eta}_2,\ldots$ of $H$.
Then $p\bar{\eta}_n$ are sections of $\E\to X$ 
and the set $\set{p(x)\bar{\eta}_n}_{n\in\N}$ is dense in $\E_x$ for every $x$. 
We will now construct bundle projections $p_1\geq p_2\geq\ldots$ of $H_X$ 
and sections $\xi_n$ of $S(p_n)$ by induction as follows.
We start with $p_1=p$. 
Suppose that $p_n$ is already constructed.
Applying the lemma above to $\eps_n=1/n$, $p_n$, and $\eta_n=p_n\bar{\eta}_n$,
we obtain a continuous section $\xi_n\colon X\to S(p_n)$ such that the distance from $\eta_n(x)$ 
to $\bra{\xi_n(x)}$ is less than $\eps_n$ for every $x\in X$.
Let $q_n(x)$ be the orthogonal projection onto $\bra{\xi_n(x)}$.
Then the map $q_n\colon X\to\Proj(H)$ is norm continuous, 
and we proceed with the next bundle projection $p_{n+1}=p_n-q_n$.
In such a way, we obtain continuous sections $\xi_1,\xi_2,\ldots$ of $S(\E)$ 
such that the vectors $\xi_n(x)$ are pairwise orthogonal 
and their finite linear combinations are dense in $\E_x$ for every $x\in X$.

The same reasoning applied to $1-p$ instead of $p$ provides continuous sections $\xi'_1,\xi'_2,\ldots$ 
of $S(H_X\ominus\E)$ such that for every $x$ the vectors $\xi'_n(x)$ 
form an orthonormal basis of $H\ominus\E_x$.

Fix an orthonormal basis $e_1,e_2,\ldots$ of $E$ and an orthonormal basis $e'_1,e'_2,\ldots$ of $H\ominus E$.
Let $u(x)\colon H\to H$ be the isometry given by the rule $u(x)e_n=\xi_n(x)$ and $u(x)e'_n=\xi'_n(x)$.
Then the range of $u(x)$ is dense in $H$, so $u(x)$ is unitary.
By the construction, $x\mapsto u(x)$ is strongly continuous and $u(x)E=\E_x$, 
so the map $u\colon X\to\sU(H)$ lifts $p\colon X\to\sP(H)$.

The trivialization of the trivial Hilbert bundle $H_X$ given by $u\inv$ takes $\E$ 
to the constant subbundle $E_X$ of $H_X$ and thus is adapted to $\E$.
\endproof

\sub{Proof of Theorem \ref{Thm:Gr-fd}.}
The structure group $\sU$ of a locally trivial Hilbert bundle $\H$ is contractible \cite[Lemma 3]{DD}.
Since the base space $X$ is paracompact, $\H$ is trivial.

Let $(X_i,\E_i)$ be an atlas of a polarization $\Pi$.
Let $x\in X_i$ be an arbitrary point.
By the assumption of the theorem, $x$ has a neighborhood $X'$ of finite covering dimension.
Since $X$ is normal, there is a closed neighborhood $Y$ of $x$ which is contained in $X'\cap X_i$.
Such an $Y$ is a closed subset of a paracompact space $X$, so it is also paracompact;
it is a closed subset of $X'$, so $\dim Y\leq\dim X'<\infty$.
By Theorem \ref{thm:Thm5-DD}, the restriction of $\E_i$ to $Y$ is trivial.
Since $x$ was chosen arbitrarily, each $\E_i$ is tame.
Therefore, $\Pi$ is also tame,
so it satisfies the assumptions of Theorem \ref{Thm:Gr} and thus admits a global adapted trivialization.
\endproof

\begin{rem}
	For a general space $X$, such a result is no longer true. 
  Indeed, Dixmier and Douady constructed a subbundle $\E$ of the trivial Hilbert bundle $\H$ 
	over a compact space $X$,	which is not locally trivial. 	
	We recall their construction in Example \ref{exm:DD} below.
	Every fiber of $\E$ has infinite rank and corank, so $\E$ determines a polarization of $\H$. 
	However, $\E$ has no adapted trivializations, even locally.
\end{rem}

\section{Self-adjoint operators}\label{sec:sa-triv}

\upskip
\sub{Proof of Theorem \ref{Thm:Beusa}.}
A regular operator has compact resolvents if and only if its bounded transform is essentially unitary.
A trivialization is adapted to a family $\A$ of regular operators if and only if 
it is adapted to the bounded transform $\f(\A)$.
Therefore, we only need to prove Theorem \ref{Thm:Beusa} for essentially unitary operators $\A_x$;
the case of operators with compact resolvents will immediately follow.

Let $X^+$, resp. $X^-$ be the subset of points $\xX$ for which $\A_x$ is essentially positive, resp. negative,
and let $X\st = X\setminus(X^+\cup X^-)$.
In a local adapted trivialization $\A$ is norm continuous, 
so each of the three subspaces $X^+$, $X^-$, and $X\st$ is open in $X$.
Therefore, $X$ is the disjoint union of $X^+$, $X^-$, and $X\st$,
so it is enough to prove the theorem for each of these subspaces separately.

For $x\in X^+$, the operator $\A_x$ is essentially positive and essentially unitary, so $\A_x-1$ is compact.
Over $X^+$ the family $x\mapsto \A_x-1$ is a norm continuous family of compact operators, 
so both $\A-1$ and $\A$ are norm continuous in every trivialization of $\H$ over $X^+$.
The similar reasoning shows that $\A$ is norm continuous in every trivialization of $\H$ over $X^-$.
Trivializations of $\H$ over $X^+\cup X^-$ are identified with maps $X^+\cup X^-\to\sU(H)$,
so every two trivializations are homotopic.

It remains to prove the theorem over $X\st$.
Without loss of generality we can suppose that $X\st=X$.
Consider the family of spectral projections
\begin{equation}\label{eq:Xla}
	p_{\lambda}(x) = \onelp (\A_x) \txt{for}
	x\in X_{\lambda} = \sett{x\in X}{\lambda\in\Res(\A_x)}. 
\end{equation}
The sets $X_{\lambda}$ with $\lambda$ running the interval $(-1,1)$ form an open covering of $X$.
A local trivialization of $\H$ adapted to $\A$ takes $p_{\lambda}$ to a norm continuous family over $X_{\lambda}$.
For $\lambda<\mu$ the difference 
\begin{equation}\label{eq:pl-pm}
	p_{\lambda}(x)-p_{\mu}(x) = \one_{[\lambda,\mu)}(\A_x) = \one_{[\lambda,\mu]}(\A_x) 
\end{equation}
is a norm continuous family of finite rank operators over $X_{\lambda}\cap X_{\mu}$.
Hence the collection $(X_{\lambda},p_{\lambda})$, $\lambda\in(-1,1)$ 
is an atlas of a tame polarization $\Pi_\A$ of $\H$ over $X$.
The difference $k_{\lambda} = \A-(2p_{\lambda}-1)$ is a norm continuous family 
of compact operators over $X_{\lambda}$ in a local trivialization adapted to $\A$,
so it is norm continuous in every trivialization of $\H$ over $X_{\lambda}$.
A trivialization of $\H$ over $X_{\lambda}$ is adapted to $\A = k_{\lambda}+(2p_{\lambda}-1)$ 
if and only if it is adapted to $p_{\lambda}$.
Hence a trivialization is adapted to $\A$ if and only if it is adapted to $\Pi_{\A}$.
It remains to apply Theorem \ref{Thm:Gr}.
\endproof

\sub{Spectral projections.}
For a family $\A_x$, $\xX$ of self-adjoint operators and $\lambda\in\R$ 
we consider the family $p_{\lambda}$ of spectral projections over the subset 
$X_{\lambda}\subset X$ defined by \eqref{eq:Xla}.

\begin{prop}\label{prop:p-g2s}
If $\A\colon X\to\Rsa(H)$ is graph continuous, 
then $p_{\lambda}\colon X_{\lambda}\to\Proj(H)$ is strongly continuous for every $\lambda\in\R$.
\end{prop}

\proof
Without loss of generality we can suppose that $X=\gR\sa(H)$ and $\A$ is the identity map.
In other words, we need to prove that the spectral projection 
\[ p_{\lambda} = \onelp\colon\Rsa_{\lambda}\to\Proj(H), \quad 
   \Rsa_{\lambda} = \sett{A\in\Rsa(H)}{\lambda\in\Res(A)} \] 
is graph-to-strong continuous.
Let $A_0\in\Rsa_{\lambda}$.
Then $\Res(A_0)$ contains the interval $[\lambda-\eps,\lambda+\eps]$ for some $\eps>0$.
By \cite[Theorem IV.3.1]{Kato}, there is a neighborhood $V$ of $A_0$ in $\gR\sa$ 
consisting of operators whose resolvent set contains this interval. 
Let $f\colon\R\to[0,1]$ be a continuous function 
equal to 0 on $(-\infty,\lambda-\eps]$ and to 1 on $[\lambda+\eps,\infty)$.
The function $g=\f-(2f-1)$ is continuous on $\R$ and vanishes at $\pm\infty$,
so it can be approximated by polynomials in $(x+i)\inv$ and $(x-i)\inv$
and thus induces a graph-to-norm continuous map $g\colon\Rsa\to\B$.
See \cite{RS}, Theorem VIII.20(a).
By Proposition \ref{prop:f-g2s} $\f$ is graph-to-strong continuous,
so $f = (\f-g+1)/2$ is also graph-to-strong continuous.
Since $p_{\lambda}$ coincides with $f$ on $V$ and $A_0$ was chosen arbitrary, 
$p_{\lambda}$ is graph-to-strong continuous on the whole $\Rsa_{\lambda}$.
\endproof

\begin{prop}\label{prop:ad-pr}
Let $\A\colon X\to\Rsa(H)$ be graph continuous and $\lambda\in\R$.
Then $\A$ is Riesz continuous on $X_{\lambda}$ if and only if 
$p_{\lambda}\colon X_{\lambda}\to\Proj(H)$ is norm continuous.
\end{prop}

\proof
Without loss of generality we can suppose that $X_{\lambda}=X$,
that is, $\lambda$ belongs to the resolvent set of each $\A_x$, $\xX$.

Suppose that $\A$ is Riesz continuous. 
Then the family $a_x=\f(\A_x)$ is norm continuous and $\mu=\f(\lambda)$ lies in the resolvent set of each $a_x$.
Therefore, the map $x\mapsto p_{\lambda}(x) = \one_{[\mu,+\infty)}(a_x)$ is also norm continuous.

Conversely, suppose that $p_{\lambda}$ is norm continuous.
Then $r_{\lambda} = 2p_{\lambda}-1$ is a norm continuous family of symmetries 
(that is, self-adjoint unitary operators) commuting with $\A$.
By Proposition \ref{prop:UR}, the family $\A'=\A\cdot r_{\lambda}$ of self-adjoint operators 
is graph (resp Riesz) continuous if and only if $\A$ is graph (resp Riesz) continuous.
But, in contrast with $\A$, the family $\A'$ is uniformly bounded from below by $-|\lambda|$.
The graph and Riesz topology coincide on the set of operators bounded from below by a given constant   
\cite[Proposition 11.2]{Pr20}, so $\A'$ and thus also $\A$ are Riesz continuous.
\endproof

\sub{Polarization associated with operator family.}
Let $\A$ be a graph continuous family of self-adjoint operators $\A_x\colon\H_x\to\H_x$ with compact resolvents. 
Then $X_{\lambda}$ is open for every $\lambda\in\R$ by \cite[Theorem IV.3.1]{Kato} 
and these sets form an open covering of $X$.
By Proposition \ref{prop:p-g2s}, $p_{\lambda}$ is a bundle projection of $\H$ over $X_{\lambda}$.
For $\lambda<\mu$ the difference \eqref{eq:pl-pm}
is a norm continuous family of finite rank operators 
over $X_{\lambda}\cap X_{\mu}$ by \cite[Theorem IV.3.16]{Kato}.
Suppose additionally that each $\A_x$ is neither essentially positive nor essentially negative.
Then the collection $(X_{\lambda},p_{\lambda})$ is an atlas of a (finite) polarization. 
Hence a family $\A$ satisfying assumptions above determines a polarization of $\H$ 
which we call the \emph{polarization associated with} $\A$ and denote by $\Pi_{\A}$.

\begin{prop}\label{prop:ad-pol}
Let $\A$ be a graph continuous family of operators $\A_x\in\Ro_K(\H_x)$ acting on fibers of a Hilbert bundle $\H$.
A trivialization of $\H$ is adapted to $\A$ if and only if 
it is adapted to the associated polarization $\Pi_{\A}$.
\end{prop}

\proof
This follows immediately from Proposition \ref{prop:ad-pr}.
\endproof

\sub{Proof of Theorem \ref{Thm:RKsa-fd}.}
Let $\Pi_{\A}$ be the polarization of $\H$ associated with $\A$.
By Proposition \ref{prop:ad-pol}, a (local or global) trivialization of $\H$ is adapted to $\A$ 
if and only if it is adapted to $\Pi_{\A}$.
It remains to apply Theorem \ref{Thm:Gr-fd}.
\endproof

\begin{rem}\label{rem:Fr-sa}
More generally, one can consider a family $\A$ of \emph{Fredholm} self-adjoint operators acting on fibers of $\H$.
Suppose that $\A$ admits local trivializations taking it to a graph continuous family.
Then the associated polarization $\Pi_{\A}$ is defined in a similar manner. 
The only change is in the definition of the subset $X_{\lambda}$:
one should add the condition that the interval $[\lambda,0]$ or $[0,\lambda]$, 
depending on the sign of $\lambda$, is disjoint from the essential spectrum of $\A_x$.

If a trivialization is adapted to such a family $\A$, then it is adapted to $\Pi_{\A}$ as well.
However, the converse statement is no longer true:
a trivialization adapted to $\Pi_{\A}$ is not necessarily adapted to $\A$,
because it might turn $\A$ into a graph discontinuous family 
and then Proposition \ref{prop:ad-pr} cannot be applied. 
This problem arises already for continuous maps $X\to\gR\st_F(H)$. 
This is the place where condition of having compact resolvents comes into the picture: 
the graph continuity of a family of operators with compact resolvents 
is independent of the choice of a trivialization.
\end{rem}

\upskip
\begin{exm}\label{exm:DD}
Let us show that Theorem \ref{Thm:RKsa-fd} does not hold without additional assumptions on the base space $X$. 
We will construct a graph continuous family $\A\colon X\to\Ro_K(H)$ 
of invertible self-adjoint operators with compact resolvents over a compact metrizable space $X$ 
such that, for every point $\xX$, $\A$ has no adapted trivialization in any neighborhood of $x$.

Our example is based on a Dixmier--Douady's construction of a non-trivial separable continuous field $\E$ 
of Hilbert spaces over a compact space $X$. See \cite[Theorem 6]{DD}.
Let us recall their construction. 
Let $H_i=\C^2$, $X_i = \PP(H_i)\cong \C\PP^1$, $H = \oplus_{i=1}^{\infty} H_i$, 
and $X = \prod_{i=1}^{\infty} X_i$ with the product topology.
The space $X$ is a countable product of metrizable spaces and thus is metrizable itself; 
it is a product of compact spaces and thus is compact itself.

For $x = (x_i)\in X$, let $\E_x\subset H$ be the closed span of $\oplus_{i=1}^{\infty} x_i$,
where $x_i\in X_i$ is considered as a one-dimensional subspace of $H_i\subset H$.
Let $p_x\in\Proj(H)$ be the projection with the range $\E_x$.
The corresponding map $p\colon X\to\Proj(H)$ is strongly (but not uniformly) continuous, 
so the subspaces $\E_x\subset H$ give rise to the subbundle (in our sense) $\E$ of $H_X$.
By \cite[Theorem 6]{DD}, every section of $\E$ vanishes at some point.
Therefore, $\E$ is not a locally trivial Hilbert bundle
(since every locally trivial Hilbert bundle over a paracompact base space is trivial 
and thus have non-vanishing sections).

Moreover, $\E$ looks the same near each point of $X$
(more precisely, the group $\prod_{i=1}^{\infty}\U(H_i)$ acts transitively on $X$
and this action lifts to the action on $\H$ by bundle automorphisms preserving $\E$),
so a local triviality of $\E$ near one point would imply local triviality over the whole $X$.
Therefore, for every open subset $X'\subset X$, $\E$ is not trivial over $X'$
and thus has no adapted trivializations over $X'$.

Let $(c_i)_{i\in\N}$ be a sequence of positive numbers converging to $\infty$.
Then 
\[ A = \oplus_{i=1}^{\infty} c_i\cdot \Id_{H_i} \] 
is a positive regular operator on $H$ with compact resolvents.
The symmetry $r_x=2p_x-1$ commutes with $A$ and 
$\A_x = A\cdot r_x$ is an invertible self-adjoint operator with compact resolvents,
which is neither essentially positive nor essentially negative.
Since $r_x$ is a strongly continuous family of unitary operators, 
the family $\A_x$ is graph continuous by Proposition \ref{prop:UR}.

The positive spectral projection of $\A_x$ is $p_x$.
By Proposition \ref{prop:ad-pr}, each local trivialization adapted to $p$ is also adapted to $\A$.
Since $p$ has no adapted trivializations,
the family $\A$ has no adapted trivializations as well.
Moreover, for every open subset $X'\subset X$, the restriction of $\A$ to $X'$ has no adapted trivialization.
\end{exm}

\begin{rem*}
	If one multiplies the family $\A_x$ from this example 
	by the strongly continuous family of unitary operators $r_x$, then the resulting family 
	$x\mapsto\A_x r_x = r_x\A_x = A$ will be Riesz continuous (and even constant).
	However, for every strongly continuous family of unitary operators $v_x\in\U(H)$, 
	the family $x\mapsto v_x\A_x v_x^*$ is Riesz discontinuous at every point $\xX$. 
	This demonstrates the difference between the self-adjoint and non-self-adjoint cases,
	that is, between one-sided multiplication and conjugation by unitary operators,
	which we discussed in the Introduction.
\end{rem*}

\end{document}